\documentclass{amsart}
\usepackage{amsmath}
\usepackage{amssymb}
\usepackage{graphicx}
\input xy
\xyoption{all}
\newtheorem{Lemma}{Lemma}[section]
\newtheorem{Th}[Lemma]{Theorem}
\newtheorem{Prop}[Lemma]{Proposition}

\newtheorem{Cor}[Lemma]{Corollary}

\newtheorem{Ex}[Lemma]{Example}

\newtheorem{Remark}[Lemma]{Remark}
\newenvironment{Proof}{{\sc Proof.}\ }{~\rule{1ex}{1ex}\vspace{0.2truecm}}
\newcommand{\Cal}[1]{{\mathcal #1}}
\newcommand{\End}{\mbox{\rm End}}

\newcommand{\ty}{{\mathbf t}}

    \newcommand{\Hom}{\operatorname{Hom}}
\newcommand{\Spec}{\operatorname{Spec}}
\newcommand{\Mod}{\operatorname{Mod-\!}}

\newcommand{\lMod}{\mbox{\rm -Mod}}

\newcommand{\soc}{\mbox{\rm soc}}

\DeclareMathOperator{\Ob}{Ob}

\newcommand{\N}{\mathbb N}
\newcommand{\Z}{\mathbb{Z}}

\newcommand{\T}{\mathbb{T}}
\newcommand{\Q}{\mathbb{Q}}
\newcommand{\R}{\mathbb{R}}

\begin{document}
    \title[Direct Products of Modules Whose Endomorphism Rings have \dots]{Direct Products of Modules Whose Endomorphism Rings have at Most Two Maximal Ideals}
        \author[Adel Alahmadi]{Adel Alahmadi}
 \address{Adel Alahmadi, Department of Mathematics, Faculty of Science,  King Abdulaziz University, Jeddah 21589, Saudi Arabia}
\email{analahmadi@kau.edu.sa}

    \author[Alberto Facchini]{Alberto Facchini}
    \address{Alberto Facchini, Dipartimento di Matematica Universit\`a di Padova, 35121 Padova, Italy}
   \email{facchini@math.unipd.it}
\thanks{Partially
supported by Universit\`a di Padova (Progetto ex 60\% ``Anelli e categorie di moduli'') and Fondazione Cassa di Risparmio di Padova e Rovigo (Progetto di Eccellenza ``Algebraic structures and their applications''.)}

      \begin{abstract} Monogeny classes and epigeny classes have proved to be useful in the study of direct sums of uniserial modules and other classes of modules. In this paper, we show that they also turn out to be useful in the study of direct products. \end{abstract}
    \maketitle

\section{Introduction}

Two right $R$-modules $M$ and $N$ are said to {\em
belong to the same monogeny class\/} (written $[M]_m=[N]_m$) if there exist a monomorphism $M\to N$ and a
monomorphism
$N\to M$.
Dually, $M$ and $N$ are said to {\em belong to
the same epigeny class\/} ($[M]_e=[N]_e$) if there exist an epimorphism $M\to N$ and
an epimorphism
$N\to M$. Recall that a module is {\em uniserial} if its lattice of submodules is linearly ordered under inclusion. In \cite[Theorem~1.9]{TAMS}, it was proved that if
$U_1$, $\dots,$
$U_n$,
$V_1$, $\dots,$ $V_t$ are non-zero uniserial right $R$-modules, then
$U_1\oplus\dots\oplus U_n\cong V_1\oplus\dots\oplus V_t$ if and
only if
$n=t$
and there are two permutations $\sigma,\tau$ of $\{1,2,\dots,n\}$ such that
$[U_i]_m=[V_{\sigma(i)}]_m$ and $[U_i]_e=[V_{\tau(i)}]_e$ for every $i=1,2,\dots,
n$. This result, which made possible the solution of a problem \cite[p.~189]{Warf} posed by Warfield in 1975, was then generalized in various directions. On the one hand, it was extended to the case of arbitrary, non-necessarily finite, families $\{\,U_i\mid i\in I\,\}$, $\{\,V_j\mid j\in J\,\}$ of uniserial modules \cite[Theorem~2.6]{Pri} (see Theorem~\ref{3.2} below). On the other hand, it was shown that similar theorems hold not only for uniserial modules, but also for cyclically presented modules over a local ring $R$, for kernels of morphisms between indecomposable injective modules, for couniformly presented modules, and more generally, for several classes of modules with at most two maximal right ideals (see \cite[Section~5]{Bull} and  \cite{AlbPav5}).

In this paper, we prove that a similar result holds not only for direct sums, but also for direct products of arbitrary families $\{\,U_i\mid i\in I\,\}$, $\{\,V_j\mid j\in J\,\}$ of uniserial modules. We show (Theorem~\ref{suff'}) that if there exist two bijections $\sigma,\tau\colon
I\to J$ such that $[U_i]_m=[V_{\sigma(i)}]_m$ and
$[U_i]_e=[V_{\tau(i)}]_e$ for every $i\in I$, then $\prod_{i\in I}U_i\cong \prod_{j \in J}V_j$. In fact, the theorem we prove is much more general, and involves completely prime ideals in categories of modules whose endomorphism rings have at most two maximal right ideals (Theorem~\ref{suff}). This allows us to apply our theorem not only to uniserial modules, but also to several other classes of modules, like the class of cyclically presented modules over a local ring and the class of kernels of morphisms between indecomposable injective modules (Section~\ref{App}).

We then show with some examples that in general  it is not possible to reverse our result (find the converse of it). It is possible to reverse it only in the particular case of slender modules (Theorems~\ref{nec} and~\ref{buupbp}). For this class of modules, it is possible to argue as in the recent paper \cite{Franetic}.

The rings we deal with are associative rings with identity $1\ne 0$, and modules are unitary modules.

\section{The main result}

In order to present our result in the most general setting, that of modules whose endomorphism rings have at most two maximal right ideals, we adopt the point of view of \cite[Section 6]{AlbPav5}. Thus, let $R$ be an associative ring with identity and $\Mod R$ the category of all right $R$-modules.
Let $\Cal C$ be a full subcategory
of $\Mod R$ whose class of objects $\Ob(\Cal C)$ consists of indecomposable right $R$-modules. Recall that a {\em completely prime ideal} $\Cal P$ of $\Cal C$ consists of a subgroup $\Cal P(A,B)$ of the additive abelian group $\Hom_R(A,B)$ for every pair of objects $A,B\in{}\Ob(\Cal C)$ such that: \newline
(1)~for every $A,B,C\in{}\Ob(\Cal C)$, every $f\colon A\to B$ and every $g\colon B\to C$, there holds $gf\in \Cal P(A,C)$ if and only if either $f\in \Cal P(A,B)$ or $g\in \Cal P(B,C)$; and \newline
(2) $\Cal P(A,A)$ is a proper subgroup of $\Hom_R(A,A)$ for every object $A\in{}\Ob(\Cal C)$.

If $A,B$ are objects of $\Cal C$, we say that $A$ and $B$ {\em belong to the same $\Cal P$ class}, and write $[A]_{\Cal P}=[B]_{\Cal P}$,  if  there exist $f\colon A\to B$ and $g\colon B\to A$ such that $f\notin\Cal P(A,B)$ and $g\notin\Cal P(B,a)$, that is, if $\Cal P(A,B)\ne\Hom(A,B)$ and $\Cal P(B,A)\ne\Hom(B,A)$. The full subcategory $\Cal C$ of $\Mod R$ is said to {\em satisfy Condition (DSP)} (direct summand property) if whenever $A,B,C,D$ are right $R$-modules with $A\oplus B\cong C\oplus D$ and $A,B,C\in{}\Ob(\Cal C)$, then also $D\in{}\Ob(\Cal C)$.

\medskip

We begin by recalling the Weak Krull-Schmidt Theorem, the  proof of which can be found in  \cite[Theorem 6.2]{AlbPav5}, followed by a preparatory lemma.

\begin{Th}\label{WKST}{\rm (Weak Krull-Schmidt Theorem)} Let $\Cal C$ be a full subcategory of $\Mod R$ in which all objects are indecomposable right $R$-modules and let $\Cal P,\Cal Q$ be two completely prime ideals of $\Cal C$ with the property that, for every $A\in{}\Ob(\Cal C)$, $f\colon A\to A$ is an automorphism if and only if $f\notin \Cal P(A,A)\cup \Cal Q(A,A)$.
Let $A_1,\dots,A_n$, $B_1,\dots,B_t$ be objects of $\Cal C$.  Then the modules $A_1 \oplus \dots \oplus A_n$ and $B_1
\oplus \dots \oplus B_t$ are isomorphic if and only
if $n=t$ and there are two permutations $\sigma, \tau$ of $\{ 1, 2,
\dots, n \}$ with $[A_i]_{\Cal P} = [B_{\sigma(i)}]_{\Cal P}$ and
$[A_i]_{\Cal Q} = [B_{\tau(i)}]_{\Cal Q}$ for every $i=1,\dots,n$.
\end{Th}

\begin{Lemma}\label{2.2} Let $\Cal C$ be a full subcategory of $\Mod R$ in which all objects are indecomposable right $R$-modules and let $\Cal P,\Cal Q$ be a pair of completely prime ideals of~$\Cal C$ with the property that, for every $A\in{}\Ob(\Cal C)$, $f\colon A\to A$ is an automorphism if and only if $f\notin \Cal P(A,A)\cup \Cal Q(A,A)$. Assume that $\Cal C$ satisfies Condition (DSP). Let $A,B,C\in{}\Ob(\Cal C)$ with $[C]_{\Cal P}=[A]_{\Cal P}$ and $[C]_{\Cal Q}=[B]_{\Cal Q}$.
Then there exists $D\in{}\Ob(\Cal C)$ with $A\oplus B\cong C\oplus D$. Moreover, $[D]_{\Cal P}=[B]_{\Cal P}$ and $[D]_{\Cal Q}=[A]_{\Cal Q}$.\end{Lemma}

\begin{Proof} Let $A,B,C$ be objects of $\Cal C$ such that $[C]_{\Cal P}=[A]_{\Cal P}$ and $[C]_{\Cal Q}=[B]_{\Cal Q}$. Then there exist morphisms $f\colon C\to A$, $g\colon A\to C$, $h\colon C\to B$ and $\ell\colon B\to C$ such that $f\notin\Cal P(C,A)$, $g\notin\Cal P(A,C)$, $h\notin\Cal Q(C,B)$ and $\ell\notin\Cal Q(B,C)$. Thus $gf\notin\Cal P(C,C)$ and $\ell h\notin\Cal Q(C,C)$.

We have four cases according to whether $gf\notin\Cal Q(C,C)$ or $gf\in\Cal Q(C,C)$ and $\ell h\notin\Cal P(C,C)$ or $\ell h\in\Cal P(C,C)$. If $gf\notin\Cal Q(C,C)$, then $gf$ is an automorphism of $C$. Since the composite mapping of $h\colon C\to B$ and $(\ell h)^{-1}\ell\colon B\to C$ is the identity mapping of $C$, it follows that $C$ is isomorphic to a direct summand of $B$. But $B$ and $C$ are indecomposable, so that $C\cong B$. In particular, $[B]_{\Cal P}=[A]_{\Cal P}$. It follows that $D:=A$ has the required properties. Similarly if $\ell h\notin\Cal P(C,C)$.

Hence it remains to consider the case $gf\in\Cal Q(C,C)$ and $\ell h\in\Cal P(C,C)$. In this case, we have that $gf+\ell h\notin\Cal P(C,C)\cup\Cal Q(C,C)$, so $gf+\ell h$ is an automorphism of $C$. Now
the composite mapping of $\binom{f}{h}\colon C\to A\oplus B$ and $(gf+\ell h)^{-1}(g\ \ell)\colon A\oplus B\to C$ is the identity mapping of $C$. Thus $C$ is isomorphic to a direct summand of $A\oplus B$. By Condition {\em (DSP)}, there exists an object $D$ of $\Cal C$ with $A\oplus B\cong C\oplus D$. Finally, $[D]_{\Cal P}=[B]_{\Cal P}$ and $[D]_{\Cal Q}=[A]_{\Cal Q}$ by Theorem~\ref{WKST}.
\end{Proof}

We are ready for the proof of the main result of this paper.

\begin{Th}\label{suff} Let $\Cal C$ be a full subcategory of $\Mod R$ in which all objects are indecomposable right $R$-modules and let $\Cal P,\Cal Q$ be two prime ideals of $\Cal C$ with the property that, for every $A\in{}\Ob(\Cal C)$, $f\colon A\to A$ is an automorphism if and only if $f\notin \Cal P(A,A)\cup \Cal Q(A,A)$. Assume that $\Cal C$ satisfies Condition (DSP).
Let $\{\,A_i\mid i\in I\,\}$ and $\{\,B_j\mid j\in J\,\}$ be two
families of objects of $\Cal C$. Assume that there exist two bijections $\sigma,\tau\colon
I\to J$ such that $[A_i]_{\Cal P} = [B_{\sigma(i)}]_{\Cal P}$ and
$[A_i]_{\Cal Q} = [B_{\tau(i)}]_{\Cal Q}$ for every $i\in I$. Then the $R$-modules $\prod_{i\in I}A_i$ and $\prod_{j \in J}B_j$ are isomorphic.\end{Th}

\begin{Proof} The proof is very similar to the proof of \cite[Theorem~3.1]{DF}. Let $\sigma,\tau\colon
I\to J$ be two bijections such that $[A_i]_{\Cal P}=[B_{\sigma(i)}]_{\Cal P}$ and
$[A_i]_{\Cal Q}=[B_{\tau(i)}]_{\Cal Q}$ for every $i\in I$.
We want to prove that $\prod_{i\in I}A_i\cong \prod_{j \in J}B_j$. Re-indexing the family $\{\,B_j\mid j\in J\,\}$ in the set $I$ via the bijection $\sigma$, we can suppose, without
 loss of generality, that $I=J$ and $\sigma\colon I\to I$ is the identity. Thus we have that $\tau$ is an element of the
symmetric group $S_I$  of all permutations of $I$, and that $[A_i]_{\Cal P}=[B_i]_{\Cal P}$ and
$[A_i]_{\Cal Q}=[B_{\tau(i)}]_{\Cal Q}$ for every $i\in I$. We must show that $\prod_{i\in I}A_i\cong \prod_{i\in I}B_i$.

The symmetric group $S_I$ acts on the set $I$.
Let $C$ be the cyclic subgroup of $S_I$ generated by $\tau$, so that $C$ acts on $I$. For every $i\in I$, let $[i]=\{\,\tau^z(i)\mid z\in \Z\,\}$
be the $C$-orbit of $i$. As $C$ is countable, the $C$-orbits are either finite or countable. We claim that $\prod_{k\in[i]}A_k\cong \prod_{k\in[i]}B_k$.
If the orbit $[i]$ is finite, in which case direct sum and direct product coincide, the claim follows immediately from Theorem~\ref{WKST}, because the modules $A_k$ and $B_k$ belong to the same $\Cal P$ class for every $k$, and the modules in the finite sets $\{\,A_k\mid k\in[i]\,\}$ and $\{\,B_k\mid k\in[i]\,\}$ have their $\Cal Q$ classes permuted by the restriction of $\tau$ to the orbit $[i]$. Hence we can assume that the orbit $[i]$ is infinite. For simplicity of notation,  set $i_z:=\tau^z(i)$, $A_z:=A_{i_z}$ and $B_z:=B_{i_z}$ for every $z\in\Z$. In this notation, we have that $\tau(i_z)=\tau^{z+1}(i)=i_{z+1}$ for every $z$, so that $[A_z]_{\Cal P}=[B_z]_{\Cal P}$ and
$[A_z]_{\Cal Q}=[B_{{z+1}}]_{\Cal Q}$ for every $z\in \Z$, and we must show that $\prod_{z\in\Z}A_z\cong \prod_{z\in\Z}B_z$.

For every integer $n\ge 0$, we will construct by induction on $n$ a 4-tuple $$(X_n,Y_n,X'_n,Y'_n)$$ of modules isomorphic to objects of $\Cal C$ with the following five properties  for every $n\ge 0$:
\begin{enumerate}
\item[{\rm (a)}] $X_n\oplus Y_n=A_{n}\oplus A_{-n-1}$;
\item[{\rm (b)}]  $X'_n\oplus Y'_n=B_{n+1}\oplus B_{-n-1}$;
\item[{\rm (c)}]  $[X_n]_{\Cal P}=[B_n]_{\Cal P}$ and $[X_n]_{\Cal Q}=[B_{-n}]_{\Cal Q}$;
\item[{\rm (d)}]  $X_{n+1}\cong X'_n$;
\item[{\rm (e)}]  $Y_n\cong Y'_n$.\end{enumerate}

As $[B_0]_{\Cal P}=[A_0]_{\Cal P}$ and $[B_0]_{\Cal Q}=[A_{-1}]_{\Cal Q}$, there is a direct-sum decomposition $X_0\oplus Y_0$ of $A_0\oplus A_{-1}$ such that $X_0\cong B_0$, $[Y_0]_{\Cal P}=[A_{-1}]_{\Cal P}$ and $[Y_0]_{\Cal Q}=[A_0]_{\Cal Q}$ (Lemma~\ref{2.2}). Hence  $[Y_0]_{\Cal P}=[B_{-1}]_{\Cal P}$ and $[Y_0]_{\Cal Q}=[B_1]_{\Cal Q}$. Thus we can apply Lemma~\ref{2.2} to the module $Y_0$ and the direct sum $B_1\oplus B_{-1}$, getting that there is a direct-sum decomposition $B_1\oplus B_{-1}= X'_0\oplus Y'_0$ of $B_1\oplus B_{-1}$ with $Y'_0\cong Y_0$ and $X'_0\in{}\Ob(\Cal C)$.
This proves that the 4-tuple $(X_0,Y_0,X'_0,Y'_0)$ has the required properties.

Now suppose $n\ge 1$ and that the 4-tuple $(X_t,Y_t,X'_t,Y'_t)$ with the required properties
has been defined for every $t$ with $0\le t<n$. Then $X_{n-1}\oplus Y_{n-1}=A_{n-1}\oplus A_{-n}$ by (a), and $[X_{n-1}]_{\Cal P}=[B_{n-1}]_{\Cal P}=[A_{n-1}]_{\Cal P}$ by (c). Thus, by Theorem~\ref{WKST}, we know that $[Y_{n-1}]_{\Cal P}=[A_{-n}]_{\Cal P}=[B_{-n}]_{\Cal P}$. Similarly, from $X_{n-1}\oplus Y_{n-1}=A_{n-1}\oplus A_{-n}$ and $[X_{n-1}]_{\Cal Q}=[B_{-n+1}]_{\Cal Q}=[A_{-n}]_{\Cal Q}$ (property (c)), we have that $[Y_{n-1}]_{\Cal Q}=[A_{n-1}]_{\Cal Q}=[B_n]_{\Cal Q}$.
Applying (e), it follows that $[Y'_{n-1}]_{\Cal P}=[B_{-n}]_{\Cal P}$ and $[Y'_{n-1}]_{\Cal Q}=[B_n]_{\Cal Q}$.
These last two equalities and $X'_{n-1}\oplus Y'_{n-1}=B_{n}\oplus B_{-n}$ imply that
\begin{equation} [X'_{n-1}]_{\Cal P}=[B_{n}]_{\Cal P}\quad\mbox{\rm and}\quad [X'_{n-1}]_{\Cal Q}=[B_{-n}]_{\Cal Q}.\label{(3)} \end{equation}
Thus $[X'_{n-1}]_{\Cal P}=[A_{n}]_{\Cal P}$ and $[X'_{n-1}]_{\Cal Q}=[A_{-n-1}]_{\Cal Q}$. By
Lemma~\ref{2.2}, there exist $X_n$ and $Y_n$ such that $X_n\oplus Y_n=A_n\oplus A_{-n-1}$  and $X'_{n-1}\cong X_n$. Thus properties (a) and (d) hold. Property (c) also follows easily.

From equalities~(\ref{(3)}), we have that $[X_n]_{\Cal P}=[X'_{n-1}]_{\Cal P}=[B_n]_{\Cal P}=[A_n]_{\Cal P}$ and $[X_n]_{\Cal Q}=[X'_{n-1}]_{\Cal Q}=[B_{-n}]_{\Cal Q}=[A_{-n-1}]_{\Cal Q}$. These equalities and $X_n\oplus Y_n=A_n\oplus A_{-n-1}$ imply that $[Y_n]_{\Cal P}=[A_{-n-1}]_{\Cal P}$ and $[Y_n]_{\Cal Q}=[A_n]_{\Cal Q}$. Hence $[Y_n]_{\Cal P}=[B_{-n-1}]_{\Cal P}$ and $[Y_n]_{\Cal Q}=[B_{n+1}]_{\Cal Q}$. These two equalities give a direct-sum decomposition $X'_n\oplus Y'_n=B_{n+1}\oplus B_{-n-1}$ with $Y_n\cong Y'_n$ (Lemma~\ref{2.2} and Condition~{\em (DSP)}). Thus properties (b) and (e) also hold, and the construction by induction is completed.

Notice that from (c) it follows that $X_0\cong B_0$.
Then $$\begin{array}{l}\prod_{z\in\Z}A_z\cong \prod_{n\ge0}A_{n}\oplus A_{-n-1}\cong \prod_{n\ge0}X_n\oplus Y_n\cong  \\ \qquad\cong X_0\oplus \prod_{n\ge0}X_{n+1}\oplus Y_n\cong B_0\oplus \prod_{n\ge0}X'_n\oplus Y'_n
\cong  \\ \qquad\cong B_0\oplus \prod_{n\ge0}B_{n+1}\oplus B_{-n-1}\cong \prod_{z\in\Z}B_z.\end{array}$$ This concludes the proof of the claim. Now the orbits $[i]$ form a partition of $I$, that is, $I$ is the disjoint union of the orbits,  and $\prod_{k\in[i]}A_k\cong \prod_{k\in[i]}B_k$ for every $i$ by the claim. Taking the direct product we conclude that $\prod_{i\in I}A_i\cong\prod_{i \in I}B_i$, as desired.
\end{Proof}

\begin{Remark} {\rm In  \cite{AdelAlbKor}, the authors have considered the condition ``the canonical functor $\Cal C\to \Cal C/\Cal P\times\Cal C/\Cal Q$ is local''. Here an additive functor $F\colon\Cal A\to\Cal B$ between preadditive categories $\Cal A$ and $\Cal B$ is said to be a {\em local} functor if, for every morphism $f\colon A\to B$ in the category $\Cal A$, $F(f)$ isomorphism in $\Cal B$ implies $f$ isomorphism in $\Cal A$.  Let us prove that if $\Cal C$ is a preadditive category, $\Cal P$ and $\Cal Q$ are two completely prime ideals of $\Cal C$, and for every $A\in{}\Ob(\Cal C)$, $f\colon A\to A$ is an automorphism if and only if $f\notin \Cal P(A,A)\cup \Cal Q(A,A)$, then the canonical functor $\Cal C\to \Cal C/\Cal P\times\Cal C/\Cal Q$ is local.

In order to see this, let $\Cal C$ be a preadditive category, $\Cal P$ and $\Cal Q$ be two completely prime ideals of $\Cal C$, and suppose that for every $A\in{}\Ob(\Cal C)$, $f\colon A\to A$ is an automorphism if and only if $f\notin \Cal P(A,A)\cup \Cal Q(A,A)$. We will apply \cite[Theorem~2.4]{AdelAlbKor} to the ideal $\Cal I:=\Cal P\cap\Cal Q$ of $\Cal C$. For any object $A\in{}\Ob(\Cal C)$, $\Cal P(A,A)\cup\Cal Q(A,A)$ is the set of all non-invertible elements of the ring $\End_R(A)$. It follows that every right (or left) ideal of $\End_R(A)$ is either contained in $\Cal P(A,A)$ or $\Cal Q(A,A)$, that is, the maximal right (or left) ideals of $\End_R(A)$ are at most $\Cal P(A,A)$ and $\Cal Q(A,A)$. In any case, $\Cal I(A,A)=\Cal P(A,A)\cap\Cal Q(A,A)$ is contained in the Jacobson radical of $\End_R(A)$. By  \cite[Theorem~2.4]{AdelAlbKor}, if $\Cal J$ is the Jacobson radical of $\Cal C$, one has that $\Cal I=\Cal P\cap\Cal Q\subseteq \Cal J$. It follows that the canonical functor $\Cal C\to \Cal C/\Cal P\times\Cal C/\Cal Q$ is local \cite[Theorem~5.3]{AdelAlbKor}.

The implication we have just proved in the previous paragraph cannot be reversed. In order to see it, consider the full subcategory $\Cal C$ of $\Mod\Z$ with the unique object $\Z$. Let $\Cal P=\Cal Q=0$ be the zero ideal of $\Cal C$, which is a completely prime ideal of $\Cal C$. The canonical functor $\Cal C\to \Cal C/\Cal P\times\Cal C/\Cal Q$ is trivially local.
Multiplication by $n\ge 2$ is an endomorphism of $\Z$ that is not an automorphism and does not belong to $0=\Cal P(\Z,\Z)\cup \Cal Q(\Z,\Z)$.}\end{Remark}

\section{Applications}\label{App}

Now we are going to apply Theorem~\ref{suff} to a number of examples.

\subsection{Biuniform modules.}\label{bm} Let $R$ be a ring and $\Cal B$ be the full subcategory of $\Mod R$ whose objects are all biuniform right $R$-modules; that is, the modules that are both uniform and couniform (=hollow; biuniform modules are those of Goldie dimension $1$ and dual Goldie dimension $1$). If $A$ and $B$ are biuniform $R$-modules, let $\Cal P(A,B)$ be the group of all non-injective morphisms $A\to B$ and $\Cal Q(A,B)$ be the group of all non-surjective morphisms $A\to B$. Then $\Cal P$ and $\Cal Q$ are completely prime ideals of $\Cal B$ \cite[Lemma~6.26]{book}, the category $\Cal B$ clearly satisfies Condition {\em (DSP)}, and the pair $\Cal P,\Cal Q$ satisfies the hypotheses of Theorems~\ref{WKST} and~\ref{suff}. Thus, from Theorem~\ref{suff}, we immediately get that:

\begin{Th}\label{suff'} Let $\{\,U_i\mid i\in I\,\}$ and $\{\,V_j\mid j\in J\,\}$ be two
families of biuniform modules over an arbitrary
ring $R$. Assume that there exist two bijections $\sigma,\tau\colon
I\to J$ such that $[U_i]_m=[V_{\sigma(i)}]_m$ and
$[U_i]_e=[V_{\tau(i)}]_e$ for every $i\in I$. Then $\prod_{i\in I}U_i\cong \prod_{j \in J}V_j$.\end{Th}

Since non-zero uniserial modules are biuniform, Theorem~\ref{suff'} holds in particular for families $\{\,U_i\mid i\in I\,\}$ and $\{\,V_j\mid j\in J\,\}$ of non-zero uniserial modules.

\subsection{Uniserial modules, quasismall modules.} Quasismall modules have a decisive role in the study of direct sums of uniserial modules. Recall that a module $N_R$ is {\em quasismall} if for every set $\{\,M_i\mid i\in I\,\}$ of $R$-modules such that
$N_R$ is isomorphic to a direct summand of $\oplus_{i\in
I}M_i$, there exists a finite subset
$F$ of $ I$ such that $N_R$ is isomorphic to a direct summand of
$\oplus_{i\in
F}M_i$. For instance, every finitely generated module is quasismall, every module with local endomorphism ring is quasismall, and every uniserial module is either quasismall or countably generated. There exist uniserial modules that are not quasismall \cite{Puni}.

Pavel Prihoda proved in \cite{Pri} (the necessity of the condition had already been proved in \cite{DF}) that:

\begin{Th}\label{3.2} Let $\{\,U_i\mid i\in I\,\}$ and $\{\,V_j\mid j\in J\,\}$ be two
families of uniserial modules over an arbitrary
ring $R$. Let $I'$ be the sets of all indices $i\in I$ with $U_i$ quasismall, and similarly for $J'$. Then $\bigoplus_{i\in I}U_i\cong \bigoplus_{j \in J}V_j$ if and only if there exist a bijection $\sigma\colon
I\to J$ such that $[U_i]_m=[V_{\sigma(i)}]_m$ and  a bijection $\tau\colon
I'\to J'$ such that
$[U_i]_e=[V_{\tau(i)}]_e$ for every $i\in I'$. \end{Th}

Therefore it is natural to ask whether Theorem~\ref{suff'} remains true for uniserial modules if we weaken its hypotheses to the condition studied by Prihoda. That is, assume that  $\{\,U_i\mid i\in I\,\}$ and $\{\,V_j\mid j\in J\,\}$ are two
families of uniserial modules over an arbitrary
ring $R$. Let $I'$ be the sets of all indices $i\in I$ with $U_i$ quasismall, and similarly for $J'$. Suppose that there exist a bijection $\sigma\colon
I\to J$ such that $[U_i]_m=[V_{\sigma(i)}]_m$ for every $i\in I$ and  a bijection $\tau\colon
I'\to J'$ such that
$[U_i]_e=[V_{\tau(i)}]_e$ for every $i\in I'$.  Is it true that $\prod_{i\in I}U_i\cong \prod_{j \in J}V_j$?

Equivalently, if $\{\,U_i\mid i\in I\,\}$ and $\{\,V_j\mid j\in J\,\}$ are two
families of uniserial modules over an arbitrary
ring $R$ and $\bigoplus_{i\in I}U_i\cong \bigoplus_{j \in J}V_j$, is  it true that $\prod_{i\in I}U_i\cong \prod_{j \in J}V_j$? We don't know what the answer to this question is, but, in a very special case, it is possible to find a result dual to Prihoda's Theorem~\ref{3.2}. For this purpose, we now recall the main results of \cite[Section~6]{FG}. Recall that if $_SA$ and $_SB$ are left modules over a ring $S$, $_SA$ is said to be {\em cogenerated by $_SB$} if $_SA$ is isomorphic to a submodule of a direct product of copies of $_SB$.

Let $R$ be any ring.  Fix a set $\{\,E_\lambda\mid\lambda\in\Lambda\,\}$
of representatives up to isomorphism of all injective right
$R$-modules that are injective envelopes of some non-zero uniserial
$R$-module.  Set $E_R:=E(\oplus_{\lambda\in\Lambda}E_\lambda)$ and
$S:=\End(E_R)$.  Then $S/J(S)$ is a von Neumann regular ring and idempotents can be lifted modulo $J(S)$, so that $S$ is an exchange ring \cite{W72}. Thus  idempotents can be lifted modulo every left (respectively right) ideal \cite{Nicho}. Moreover, ${}_SE_R$ turns out to be an $S$-$R$-bimodule and
$$\Hom(-,{}_SE_R)\colon \Mod R\to S\lMod$$ is a contravariant
exact functor. Let $\Cal C_R$ denote the full subcategory of $\Mod R$
whose objects are all uniserial right $R$-modules. Let ${}_S\Cal C'$ be the full subcategory of $S \lMod$
whose objects are all uniserial left $S$-modules
that have a projective cover and are cogenerated by $_SE$. It is possible to prove that if a non-zero uniserial module $U$ has a projective cover $P$, then $P$ is a couniform module \cite[Lemma~2.2]{FG}, so that, in particular, $P$, hence $U$, are cyclic modules. Thus all the $S$-modules in ${}_S\Cal C'$ are quasismall. The following result is proved in \cite[Proposition~6.1 and Corollary~6.2]{FG}.

\begin{Prop}\label{vhl}
{\rm (a)} The restriction $H$ of the functor $\Hom(-,{}_SE_R)\colon \Mod R\to S\lMod$ is a duality between the categories $\Cal C_R$ and ${}_S\Cal C'$.

{\rm (b)} Two uniserial right $R$-modules $U_R$ and $U'_R$ belong to the same monogeny class if and only if the uniserial left $S$-modules  $H (U_R)$ and $H(U'_R)$ belong to the same epigeny class.

{\rm (c)} Two uniserial right $R$-modules   $U_R$ and $U'_R$ belong to the same epigeny class if and only if the uniserial left $S$-modules  $H(U_R)$ and $H(U'_R)$ belong to the same monogeny class.\end{Prop}

Call {\em dually quasismall} any object of ${}_S\Cal C'$ isomorphic to $H(U_R)$ for some quasismall uniserial $R$-module $U_R$. Then we have that:

\begin{Th}\label{Pridual} Let
$\{\,U_i\mid i\in I\,\}$, $\{\,V_j\mid j\in J\,\}$ be two
families of objects of ${}_S\Cal C'$. Let $I'$ be the sets of all the indices $i\in I$ with $U_i$ dually quasismall, and $J'$ be the sets of all the indices $j\in J$ with $V_j$ dually quasismall. Suppose that there exist a bijection $\sigma\colon
I'\to J'$ such that $[U_i]_m=[V_{\sigma(i)}]_m$ for every $i\in I'$ and  a bijection $\tau\colon
I\to J$ such that
$[U_i]_e=[V_{\tau(i)}]_e$ for every $i\in I$.  Then $\prod_{i\in I}U_i\cong \prod_{j \in J}V_j$.\end{Th}

\begin{Proof} For every $i\in I$, let $X_i$ be a uniserial right $R$-module with $H(X_i)\cong U_i$ and, for every $j\in J$, let $Y_j$ be a uniserial right $R$-module with $H(Y_j)\cong V_j$. The module $X_i$ is quasismall if and only if $i\in I'$, and $Y_j$ is quasismall if and only if $j\in J'$. By Proposition~\ref{vhl}, $[X_i]_e=[Y_{\sigma(i)}]_e$ for every $i\in I'$ and
$[X_i]_m=[Y_{\tau(i)}]_m$ for every $i\in I$.  From Theorem~\ref{3.2}, we know that
$\bigoplus_{i\in I}X_i\cong \bigoplus_{j \in J}Y_j$. Applying the functor $\Hom(-,{}_SE_R)$, we obtain the desired conclusion $\prod_{i\in I}U_i\cong \prod_{j \in J}V_j$.\end{Proof}

\subsection{Cyclically presented modules}

For any ring $R$, we denote by $U(R)$ the group of all invertible elements of $R$ and by $J(R)$ the Jacobson radical of $R$.
Let $R$ be a local ring and $\Cal C$ be the full subcategory of $\Mod R$ whose objects all the modules $R/rR$ with $r\in J(R)\setminus\{0\}$. Since $R$ is local, all the modules in $\Cal C$ are couniform, and therefore all objects of $\Cal C$ are indecomposable modules. If $R/rR,R/sR\in\Ob(\Cal C)$, every morphism $R/rR\to R/sR$ is induced by left multiplication by some element $t\in R$, and $\Hom(R/rR,R/sR)\cong \{\, t\in R\mid tr\in sR\,\}/sR$.
Let $\Cal P(R/rR,R/sR)$ be the group $\{\, t\in R\mid tr\in sJ(R)\,\}/sR$. Then $\Cal P$ turns out to be a completely prime ideal of $\Cal C$.
If $R/rR$ and $R/sR$ are objects of $\Cal C$ in the same $\Cal P$ class, the modules $R/rR$ and $R/sR$ are said to  {\em have the same lower part\/}, denoted by $[R/rR]_l=[R/sR]_l$ \cite{Amini}. It is easily seen that $[R/rR]_l=[R/sR]_l$ if and only if there exist $u,v\in U(R)$ and $x,y\in R$ with $ru= xs$ and $sv= yr$. As in the previous Example~\ref{bm}, let $\Cal Q(R/rR,R/sR)$ be the group of all non-surjective morphisms $R/rR\to R/sR$. In this case, $\Cal Q(R/rR,R/sR)\cong \{\, t\in J(R)\mid tr\in sR\,\}/sR$.

The pair $\Cal P,\Cal Q$ of completely prime ideals satisfies the hypotheses of  Theorems~\ref{WKST} and~\ref{suff}, and the category $\Cal C$ satisfies Condition {\em (DSP)}, so that, from Theorem~\ref{suff}, we immediately obtain:

\begin{Th}\label{suff''} Let $R$ be a local ring and $\{\,U_i\mid i\in I\,\}$ and $\{\,V_j\mid j\in J\,\}$ be two
families of right $R$-modules in $\Cal C$. Suppose that there exist two bijections $\sigma,\tau\colon
I\to J$ such that $[U_i]_l=[V_{\sigma(i)}]_l$ and  and
$[U_i]_e=[V_{\tau(i)}]_e$ for every $i\in I$. Then $\prod_{i\in I}U_i\cong \prod_{j \in J}V_j$.\end{Th}

More generally, lower part and epigeny class can be defined for couniformly presented right modules over a (non-necessarily local) ring. Here an $R$-module $M$ is said to be {\em couniformly presented} if there exists an exact sequence $0\to M_1\to P\to M\to 0$ with $P$ projective and $P$ and $M_1$ of dual Goldie dimension $1$; cf.~\cite{FG}. Also in this case, we have an analogue of Theorem~\ref{suff''} for couniformly presented $R$-modules.

\subsection{Kernels of morphisms between indecomposable injective modules.}

Let $A$ and $B$ be two modules. We say that $A$ and $B$ {\it have the same upper part},
and write $[A]_u=[B]_u$, if there exist a homomorphism $\varphi\colon
E(A)\rightarrow E(B)$ and a homomorphism $\psi\colon E(B)\rightarrow E(A)$
such that $\varphi^{-1}(B)=A$ and $\psi^{-1}(A)=B$.
If $E_1,E_2,E'_1,E'_2$ are injective indecomposable right modules over an arbitrary ring $R$ and $\varphi\colon E_1\to E_2$, $\varphi'\colon E'_1\to E'_2$ are arbitrary morphisms, then $\ker\varphi \cong \ker\varphi '$ if and only if\/ $[\ker\varphi ]_m=[\ker\varphi ']_m$ and\/ $[\ker\varphi ]_u=[\ker\varphi ']_u$ \cite[Lemma~2.4]{Tufan}.

Let $R$ be an arbitrary ring and $\Cal K$ be the full subcategory of $\Mod R$ whose objects are all kernels of morphisms $f\colon E_1\to E_2$, where $E_1$ and $E_2$ range in the class of all uniform injective modules. The canonical functor $P\colon \Mod R\to\Spec(\Mod R)$, where $\Spec(\Mod R)$ denotes the spectral category of $\Mod R$ \cite{gabrieloberst}, is a left exact, covariant, additive functor, which has an $n$-th right derived functor $P^{(n)}$ for every $n\ge 0$ \cite[Proposition~2.2]{Chennai}. The restriction of $P^{(1)}$ to $\Cal K$ is a functor $\Cal K\to \Spec(\Mod R)$. If $\Cal Q(A,B)$ consists of all  morphisms $f\colon A\to B$ in $\Cal K$ with $P^{(1)}(f)=0$, then $\Cal Q$ is a completely prime ideal of $\Cal K$. If $\Cal P$ is the ideal of all non-injective homomorphisms, as in~\S\ref{bm} for biuniform modules, then the pair $\Cal P,\Cal Q$ satisfies the hypotheses of Theorems~\ref{WKST} and~\ref{suff}. In particular, the class $\Cal K$ satisfies Condition {\em (DSP)}.

Notice that in \cite{Ece} it was proved that if $\{\,A_i\mid i\in I\,\}$ and $ \{\,B_j\mid j\in J\,\}$ are two families of modules over a ring~$R$, all the $B_j$'s are kernels of non-injective morphisms between indecomposable injective modules and there exist bijections
$\sigma,\tau\colon I\to J$ such that
$[A_i]_m=[B_{\sigma(i)}]_m$ and $[A_i]_u=[B_{\tau(i)}]_u$ for every $i\in I$, then $\oplus_{i\in I}A_i\cong
\oplus_{j \in J}B_j.$ From Theorem~\ref{suff} we see that, under the same hypotheses, $\prod_{i\in I}A_i\cong
\prod_{j \in J}B_j.$

\bigskip

As a final example for this section, we can consider the following category $\Cal C$. Let $R$ be a ring and let $S_1,S_2$ be two fixed non-isomorphic simple right $R$-modules. Let $\Cal C$ be the full subcategory of $\Mod R$ whose objects are all artinian right $R$-modules $A_R$ with $\soc(A_R)\cong S_1\oplus S_2$. Set $\Cal P_i(A,B):=\{\,f\in\Hom_R(A,B)\mid f(\soc_{S_i}(A))=0\,\}$ \cite[Example~6.3(7)]{AlbPav5}. The pair of completely prime ideals $\Cal P_1,\Cal P_2$ satisfies the hypotheses of Theorem~\ref{suff}.

\section{Reversing the main result}

In this section, we consider the problem of reversng the implications in Theorems~\ref{suff} and~\ref{3.2}, that is, whether a direct product of uniserial modules determines the monogeny classes and the epigeny classes of the factors. We give four examples that prove that the answer is negative in general.

\begin{Ex}{\rm The following example shows on the one hand that it is impossible to reverse the implication in Theorem~\ref{3.2} and, on the other hand, that it is impossible to prove a result for direct products analogous to the result proved by Prihoda for direct sums of uniserial modules. More precisely, the example shows that there are two families $\{\,U_i\mid i\in I\,\}$ and $\{\,V_j\mid j\in J\,\}$ of non-zero uniserial quasismall modules over a
ring $R$ with no bijection $\sigma\colon
I\to J$ such that $[U_i]_m=[V_{\sigma(i)}]_m$  for every $i\in I$ and no a bijection $\tau\colon
I\to J$ such that
$[U_i]_e=[V_{\tau(i)}]_e$ for every $i\in I$, but with $\prod_{i\in I}U_i\cong \prod_{j \in J}V_j$. (Equivalently, $\bigoplus_{i\in I}U_i\not\cong \bigoplus_{j \in J}V_j$, but $\prod_{i\in I}U_i\cong \prod_{j \in J}V_j$.)

In this example, $R$ is the localization of the ring $\Z$ of integers at a maximal ideal $(p)$, $I=\N$ (the set of non-negative integers), $J=\N^*=\N\setminus\{0\}$ (the set of positive integers), $U_0$ is the field of fractions $\Q$ of $R$ and $U_n=V_n=\Z(p^\infty)=\Q/R$ (the Pr\"ufer group) for every $n\ge 1$. Both the $R$-modules $\Q$ and $\Z(p^\infty)$ are uniserial and with a local endomorphism ring.  Hence they are quasismall. Since $[\Q]_m\ne [\Z(p^\infty)]_m$ and $[Q]_e\ne [\Z(p^\infty)]_e$, there are no bijections $\sigma,\tau\colon \N\to\N^*$ preserving the monogeny classes and the epigeny classes, respectively. In order to show that the $R$-modules $\Q\oplus(\Z(p^\infty))^{\N^*}$ and $(\Z(p^\infty))^{\N^*}$ are isomorphic, it suffices to prove that these two divisible groups are isomorphic. Recall that two divisible abelian groups are isomorphic if and only if they have the same torsion-free rank and the same $p$-rank for every prime $p$. Thus it is enough to show that the torsion-free rank of the divisible abelian group $(\Z(p^\infty))^{\N^*}$ is infinite, that is, that the group $(\Z(p^\infty))^{\N^*}$ contains a free abelian subgroup of infinite rank. Now $$(\Z(p^\infty))^{\N^*}=(\Q/R)^{\N^*}=\{\,(q_n+R)_{n\ge 1}\mid q_n\in\Q\,\}.$$ Consider the infinitely many elements $$\left(\frac{1}{p^{nt}}+R\right)_{n\ge 1}$$ of $(\Z(p^\infty))^{\N^*}$, with $t\in\N^*$. It is easy to see that these countably many elements form a free set of generators of a free abelian subgroup of $(\Z(p^\infty))^{\N^*}$. Thus the divisible abelian group $(\Z(p^\infty))^{\N^*}$ has infinite torsion-free rank; hence the $R$-modules
$\Q\oplus(\Z(p^\infty))^{\N^*}$ and $(\Z(p^\infty))^{\N^*}$ are isomorphic. 
}\end{Ex}

\begin{Ex}{\rm Here is another example that proves that it is impossible to reverse the implication in Theorem~\ref{suff}. Let $R$ be a ring and $\Cal C$ be the full subcategory of $\Mod R$ whose objects are all injective indecomposable $R$-modules. If $A$ and $B$ are objects of $\Cal C$, let $\Cal P(A,B)$ be the group of all non-injective morphisms $A\to B$, so that $\Cal P$ is a completely prime ideal of $\Cal C$, the category $\Cal C$ satisfies Condition {\em (DSP)}, and the ideals $\Cal P=\Cal Q$ satisfy the hypotheses of Theorems~\ref{WKST} and~\ref{suff}.
Now take, for instance, $R=\Z$ and consider the family consisting of all the Pr\"ufer groups $\Z(p^\infty)$, where $p$ ranges in the set of prime numbers, and the group $\Q$. The groups in this family are pair-wise non-isomorphic and have distinct $\Cal P$ classes. In order to show that Theorem~\ref{suff} cannot be reversed, it suffices to prove that $\Q\oplus\prod_p\Z(p^\infty)\cong \prod_p\Z(p^\infty)$, because there does not exist a bijection $\sigma$ preserving the $\Cal P$ classes. And, as in the previous example, to prove that $\Q\oplus\prod_p\Z(p^\infty)\cong \prod_p\Z(p^\infty)$, it is enough to show that the divisible abelian group $\prod_p\Z(p^\infty)$ has
infinite torsion-free rank. Now the torsion subgroup of $\prod_p\Z(p^\infty)$ is the countable divisible group $\bigoplus_p\Z(p^\infty)$, so that $\prod_p\Z(p^\infty)/\bigoplus_p\Z(p^\infty)$ is a torsion-free divisible abelian group of cardinality $2^{\aleph_0}$. Thus $\prod_p\Z(p^\infty)/\bigoplus_p\Z(p^\infty)$ is a divisible group of torsion-free rank $2^{\aleph_0}$. Hence the divisible group $\prod_p\Z(p^\infty)$ has torsion-free rank $2^{\aleph_0}$. Therefore $\Q\oplus\prod_p\Z(p^\infty)\cong \prod_p\Z(p^\infty)$, though the $\Cal P$ class of $\Q$ does not appear in the set of the $\Cal P$ classes of the $\Z(p^\infty)$'s. Notice that $\prod_p\Z(p^\infty)\cong\T:=\R/\Z$.

This argument can be extended to any non-semilocal commutative Dedekind domain $R$ with cardinality $|R|=\alpha$, with maximal spectrum (set of maximal ideals) of cardinality $\beta$ for which $\alpha<\alpha^\beta$. Let $R$ be such a ring. As finite domains are fields, hence semilocal rings, it follows that $\alpha\ge\aleph_0$. Let $Q$ be the field of fractions of $R$, which also must have cardinality $\alpha$. Since all ideals in a Dedekind domain can be generated with two elements, we get that $\aleph_0\le\beta\le\alpha$. In a Dedekind domain, divisible modules coincide with injective modules. For every maximal ideal $P$ of $R$, we have that $|E(R/P)|=|Q/R_P|=\alpha$, where $R_P$ denotes the localization of $R$ at $P$. Thus $|\prod_P Q/R_P|=\alpha^\beta$. Let us prove that the torsion submodule of $\prod_P Q/R_P$ is $\bigoplus_PQ/R_P$. Let $$(q^{(P)}+R_P)_P\in \prod_P Q/R_P$$ be a torsion element. Then there exists a non-zero $r\in R$ such that $rq^{(P)}\in R_P$ for every maximal ideal $P$. Thus $Rrq^{(P)}\subseteq R_P$ for every $P$. The non-zero ideal $Rr$ of $R$ is contained in only finitely many maximal ideals of $R$, so that $R_Pr=R_P$ for almost all maximal ideals $P$. Thus $Rrq^{(P)}\subseteq R_P$ implies that $R_Prq^{(P)}\subseteq R_P$, so $R_Pq^{(P)}\subseteq R_P$, that is, $q^{(P)}\in R_P$ for almost all $P$. This proves that the torsion submodule of $\prod_P Q/R_P$ is $\bigoplus_PQ/R_P$, which
has cardinality $\alpha<\alpha^\beta$. Therefore it is possible to argue as in the previous paragraph.}\end{Ex}

\begin{Ex}{\rm This example is taken from \cite[Example 2.1]{Franetic}. Let $p$ be a prime number and $\widehat{\Z_p}$ be the ring of $p$-adic integers, so that $\Z/p^n\Z$ is a module over $\widehat{\Z_p}$ for every integer $n\ge 1$. Let $$\varphi\colon\prod_{n\ge1}\Z/p^n\Z\to \prod_{n\ge1}\Z/p^n\Z$$ be the $\widehat{\Z_p}$-module morphism defined by $\varphi(a_n+n\Z)_{n\ge1}=(a_{n+1}-a_n+n\Z)_{n\ge1}$. This morphism $\varphi$ is onto and its kernel is isomorphic to $\widehat{\Z_p}$. Thus there is an exact sequence $\xymatrix{
0 \ar[r] &  \widehat{\Z_p}\ar[r] & \prod_{n\ge1}\Z/p^n\Z\ar[r]^{\varphi} & \prod_{n\ge1}\Z/p^n\Z\ar[r] & 0,
}
$ which is a pure-exact sequence, and $\widehat{\Z_p}$ is pure-injective, so that the pure-exact sequence splits. Thus $ \widehat{\Z_p}\oplus\prod_{n\ge1}\Z/p^n\Z\cong \prod_{n\ge1}\Z/p^n\Z$. In these direct products, all factors $\widehat{\Z_p}$ and $\Z/p^n\Z$ ($n\ge1$) are pair-wise non-isomorphic uniserial $\widehat{\Z_p}$-modules, have distinct mono\-geny classes and distinct epigeny classes. Hence there cannot be bijections $\sigma$ and $\tau$ preserving the monogeny and the epigeny classes in the two direct-product decompositions. Notice that all factors have  a local endomorphism ring; hence they are quasismall uniserial modules.
}\end{Ex}

\begin{Ex}{\rm Here is a further example that shows that the condition in Theorem~\ref{suff'} is sufficient but not necessary for the isomorphism $\prod_{i\in I}U_i\cong \prod_{j \in J}V_j$ to hold. Let $U_0$ be a uniserial non-quasismall $R$-module (it is known that such modules exist  \cite{Puni}). Then $U_0$ is countably generated \cite[Lemma~4.2]{DF}; hence it is a union of an ascending chain $U_n$, $n\ge 1$ of cyclic submodules. Then $\oplus_{n\ge0}U_n\cong\oplus_{n\ge1}U_n$ \cite[Theorem~4.9]{DF}. Since $U_n$ is cyclic for $n\ge 1$ but not for $n=0$, it follows that $[U_0]_e\ne [U_n]_e$ for every $n\ge 1$. Hence there does not exist a bijection between the epigeny classes of $\{\,U_n\mid n\ge 0\,\}$ and the epigeny classes of $\{\,U_n\mid n\ge 1\,\}$. Anyway, there does exist a bijection between the monogeny classes of $\{\,U_n\mid n\ge 0\,\}$ and the monogeny classes of $\{\,U_n\mid n\ge 1\,\}$ (Theorem~\ref{3.2}). Applying the duality $H$ of Proposition~\ref{vhl}, we obtain two isomorphic direct-product decompositions $\prod _{n\ge0}H(U_n)\cong\prod _{n\ge1}H(U_n)$, with all the modules $H(U_n)$, $n\ge0$, cyclic uniserial left $S$-modules, for which there exists a bijection between the epigeny classes of $\{\,H(U_n)\mid n\ge 0\,\}$ and the epigeny classes of $\{\,H(U_n)\mid n\ge 1\,\}$, but there does not exist a bijection between the monogeny classes of $\{\,H(U_n)\mid n\ge 0\,\}$ and the monogeny classes of $\{\,H(U_n)\mid n\ge 1\,\}$. Hence the condition in Theorem~\ref{suff'} is not necessary. Notice that by Theorem~\ref{Pridual}, there is a bijection between the monogeny classes of the dually quasismall modules, that is, the modules $H(U_n)$ with $n\ge1$.}\end{Ex}

\section{Slender modules.}

Now we adopt the point of view of \cite{Franetic}, restricting our attention to slender modules. Let $R$ be a ring and $R^\omega=\prod_{n<\omega}e_nR$ be the right $R$-module that is the direct product of countably many copies of the right $R$-module $R_R$, where $e_n$ is the element of $R^\omega$ with support $\{n\}$ and equal to $1$ in $n$. A right $R$-module $M_R$ is {\em slender} if, for every homomorphism $f\colon R^\omega\to M$ there exists $n_0<\omega$ such that $f(e_n)=0$ for all $n\ge n_0$. The most important property of slender modules we need in the sequel is the following \cite[Theorem~1.2]{EM}: A module $M_R$ is slender if and only if for every countable family $\{\,P_n\mid n\ge 0\,\}$ of right $R$-modules and any homomorphism $f\colon \prod_{n\ge 0}P_n\to M_R$ there exists $m\ge0$ such that $f(\prod_{n\ge m}P_n)=0$. Here $\prod_{n\ge m}P_n$ is the subgroup of $\prod_{n\ge 0}P_n$ consisting of all elements with support contained in $\{m,m+1,m+2,\dots\}$. In the following, the cardinality of any set $I$ is denoted by $|I|$.
If $M_R$ is slender and $\{\,P_i\mid i\in I\,\}$ is a family of right $R$-modules with $|I|$ non-measurable, then $\Hom(\prod_{i\in I}P_i,M_R)\cong\bigoplus_{i\in I}\Hom(P_i,M_R)$.

Every submodule of a slender module is a slender module \cite[Lemma~1.6(i)]{EM}, so that:

\begin{Lemma} If $U_R$ is a slender module, then every module in the same monogeny class as $U_R$ is slender.\end{Lemma}

As far as Condition {\em (DSP)} is concerned, it is easily seen that:

\begin{Lemma} If $U_R, V_R$ are slender modules, then every direct summand of $U_R\oplus V_R$ is slender.\end{Lemma}

\begin{Th}\label{nec} Let $\Cal C$ be a full subcategory of $\Mod R$ in which all objects are indecomposable slender right $R$-modules and let $\Cal P,\Cal Q$ be a pair of completely prime ideals of $\Cal C$ with the property that, for every $A\in{}\Ob(\Cal C)$, $f\colon A\to A$ is an automorphism if and only if $f\notin \Cal P(A,A)\cup \Cal Q(A,A)$. Assume that $\Cal C$ satisfies Condition (DSP).
Let $\{\,A_i\mid i\in I\,\}$ and $\{\,B_j\mid j\in J\,\}$ be two
families of objects of $\Cal C$ with $|I|$ and $|J|$ non-measurable. Assume that:

{\rm (a)} In both families, there are at most countably many modules in each $\Cal P$ class.

{\rm (b)} In both families, there are at most countably many modules in each $\Cal Q$ class.

{\rm (c)} The $R$-modules $\prod_{i\in I}A_i$ and $\prod_{j \in J}B_j$ are isomorphic.

Then there exist two bijections $\sigma,\tau\colon
I\to J$ such that $[A_i]_{\Cal P} = [B_{\sigma(i)}]_{\Cal P}$ and
$[A_i]_{\Cal Q} = [B_{\tau(i)}]_{\Cal Q}$ for every $i\in I$. \end{Th}

\begin{Proof} {\em Step 1. Assume that a slender $R$-module $B$ is isomorphic to a direct summand of the direct product $\prod_{i\in I}A_i$, with $|I|$ non-measurable. Then there is a finite subset $F$ of $I$ such that $B$ is isomorphic to a direct summand of $\bigoplus_{i\in F}A_i$}.

This is \cite[Lemma 1.1]{Franetic}.

\smallskip

 {\em Step 2. For every $j\in J$ there exist $i,k\in I$ such that $[A_i]_{\Cal P}=[B_j]_{\Cal P}$ and $[A_k]_{\Cal Q}=[B_j]_{\Cal Q}$.}

Fix $j\in J$. Since $\prod_{i\in I}A_i\cong\prod_{j \in J}B_j$, we know that $B_j$ is isomorphic to a direct summand of $\prod_{i\in I}A_i$. By Step 1, there exists a finite subset $F=\{i_1,\dots,i_t\}$ of $I$ such that $B_j$ is isomorphic to a direct summand of $\bigoplus_{i\in F}A_{i}=\bigoplus_{\ell=1}^tA_{i_\ell}$. Thus there are morphisms $\varphi\colon B_j\to \bigoplus_{\ell=1}^tA_{i_\ell}$ and $\psi\colon \bigoplus_{\ell=1}^tA_{i_\ell}\to B_j$ with $\psi\varphi=1_{B_j}$. In matrix notation, we have that
$$\varphi=\left(\begin{array}{c}\varphi_1\\ \vdots \\ \varphi_t \end{array}\right)\qquad\mbox{\rm and}\qquad\psi=\left(\begin{array}{ccc}\psi_1, & \dots & ,\psi_t\end{array}\right)$$
for suitable morphisms  $\varphi_\ell\colon B_j\to A_{i_\ell}$ and $\psi_\ell\colon A_{i_\ell}\to B_j$. Thus $\sum_{\ell=1}^t\psi_\ell\varphi_\ell=1_{B_j}$. In particular, $\sum_{\ell=1}^t\psi_\ell\varphi_\ell\notin\Cal P(B_j,B_j)$. Since $\Cal P(B_j,B_j)$ is an ideal of the ring $\End_R(B_j)$, it follows that there exists an index $\ell=1,\dots,t$ such that $\psi_\ell\varphi_\ell\notin\Cal P(B_j,B_j)$. Thus $\psi_\ell\notin\Cal P(A_{i_\ell},B_j)$ and $\varphi_\ell\notin\Cal P(B_j,A_{i_\ell})$. Hence $[A_{i_\ell}]_{\Cal P}=[B_j]_{\Cal P}$, and the index $i:=i_\ell\in I$ has the required property. By symmetry, the same holds for $\Cal Q$ classes.

\smallskip

 {\em Step 3. Proof of the statement of the Theorem.}

 For every slender right $R$-module $D$, set $I_D:=\{\,i\in I\mid [A_i]_{\Cal P}=[D]_{\Cal P}\,\}$ and $J_D:=\{\,j\in J\mid [B_j]_{\Cal P}=[D]_{\Cal P}\,\}$. It suffices to show that $|I_D|=|J_D|$ for every such slender right module $D$. By contradiction, suppose that $|I_D|\ne|J_D|$ for some slender module $D$. By symmetry, we can suppose without loss of generality that $|I_D|<|J_D|$. By (a), the cardinal $|I_D|$ must be finite. Also, $J_D$ must be non-empty, so that there exists $j\in J_D$. By Step 2,
 there exist $i,k\in I$ such that $[A_i]_{\Cal P}=[B_j]_{\Cal P}$ and $[A_k]_{\Cal Q}=[B_j]_{\Cal Q}$. If $i=k$, then $A_i\cong B_j$, and since $A_i\cong B_j$ has a semilocal endomorphism ring, it can be canceled from direct sums \cite[Corollary~4.6]{book}, so that $\prod_{i'\in I}A_{i'}\cong\prod_{j" \in J}B_{j'}$ implies $\prod_{i'\in I\setminus\{i\}}A_{i'}\cong\prod_{j" \in J\setminus\{j\}}B_{j'}$. If $i\ne k$, then $B_j\oplus X_j\cong A_i\oplus A_k$ for some object $X_j$ of $\Cal C$ (Lemma~\ref{2.2}). Then $\prod_{i'\in I}A_{i'}\cong \prod_{j' \in J}B_{j'}$ can be rewritten as $B_j\oplus X_j\oplus\prod_{i'\in I\setminus\{i,k\}}A_{i'}\cong \prod_{j' \in J}B_{j'}$, and canceling as before we get that $X_j\oplus\prod_{i'\in I\setminus\{i,k\}}A_{i'}\cong \prod_{j' \in J\setminus\{j\}}B_{j'}$. In both cases $i=k$ and $i\ne k$, we have obtained two direct-product decompositions in which the families $I_D$ and $J_D$ of $\Cal P$ classes have one element less. We now proceed recursively, after $|I_D|$ steps we obtain two direct-product decompositions in which the family $I_D$ is empty and the family $J_D$ is not empty. This contradicts Step 2, and the contradiction proves that the bijection $\sigma$ with the required property exists. It is a similar argument for $\tau$.
\end{Proof}

\begin{Cor}\label{nec'} Let $\Cal C$ be a full subcategory of $\Mod R$ in which all objects are indecomposable slender right $R$-modules and let $\Cal P,\Cal Q$ be a pair of completely prime ideals of $\Cal C$ with the property that, for every $A\in{}\Ob(\Cal C)$, $f\colon A\to A$ is an automorphism if and only if $f\notin \Cal P(A,A)\cup \Cal Q(A,A)$. Assume that $\Cal C$ satisfies Condition (DSP).
Let $\{\,A_i\mid i\in I\,\}$ and $\{\,B_j\mid j\in J\,\}$ be two countable
families of objects of $\Cal C$. Assume that  $\prod_{i\in I}A_i\cong\prod_{j \in J}B_j$.Then there exist two bijections $\sigma,\tau\colon
I\to J$ such that $[A_i]_{\Cal P} = [B_{\sigma(i)}]_{\Cal P}$ and
$[A_i]_{\Cal Q} = [B_{\tau(i)}]_{\Cal Q}$ for every $i\in I$. \end{Cor}

If $\Cal D$ is a preadditive category, $A$ is an object of $\Cal D$ and $I$ is an ideal of the ring $\End_{\Cal D}(A)$, let $\Cal A_I$ be the ideal of the category $\Cal D$ defined in the following way. A morphism $f \colon X \to Y$ in $\Cal D$ belongs to $\Cal A_I(X,Y)$ if and only if $\beta f \alpha \in I$ for every pair of morphisms $\alpha \colon
A \to X$ and $\beta \colon Y \to A$ in the category $\Cal D$. The ideal $\Cal A_I$ is called the {\em ideal of $\Cal D$ associated to}~$I$ \cite{AlbPav3}. The ideal $\Cal A_I$ is the greatest of the ideals $\Cal Q$ of $\Cal D$ with $\Cal Q(A,A)\subseteq I$. It is easily seen that $\Cal A_I(A,A) = I$.

\begin{Th}\label{buupbp} Let\/ $\Cal C$ be a full subcategory of $\Mod R$
in which all objects are slender right $R$-modules
and let $\Cal P$ be a completely prime ideal of $\Cal C$. Let $\{\,A_i\mid i\in I\,\}$ and $\{\,B_j\mid j\in J\,\}$ be two
families of objects of $\Cal C$ with $|I|$ and $|J|$ non-measurable. Assume that:

{\rm (a)} For every object $A$ of $\Cal C$, $\Cal P(A,A)$ is a maximal right ideal of $\End_R(A)$.

{\rm (b)} There are at most countably many modules in each $\Cal P$ class in both families $\{\,A_i\mid i\in I\,\}$ and $\{\,B_j\mid j\in J\,\}$.

{\rm (c)} The $R$-modules $\prod_{i\in I}A_i$ and $\prod_{j \in J}B_j$ are isomorphic.

Then there is a bijection $\sigma_{\Cal P}\colon
I\to J$ such that $[A_i]_{\Cal P} = [B_{\sigma_{\Cal P}(i)}]_{\Cal P}$ for every $i\in I$. \end{Th}

\begin{Proof} Fix an object $D$ in $\Cal C$. Let
$\Cal A_{\Cal P(D,D)}$ be the ideal of $\Mod R$ associated to the maximal right ideal, hence maximal two-sided ideal, $\Cal P(D,D)$ of $\End_R(D)$. Let $P\colon\Mod R\to\Mod R/\Cal A_{\Cal P(D,D)}$ be the canonical functor. The ideal $\Cal A_{\Cal P(D,D)}$ restricts to an ideal of the category $\Cal C$, and we will also denote this restriction by $\Cal A_{\Cal P(D,D)}$. Similarly, the restriction of $P$ will be still denoted by $P$, so that $P\colon\Cal C\to\Cal C/\Cal A_{\Cal P(D,D)}$.

\smallskip

{\em Step 1. $D$ is a non-zero object in the factor category $\Cal C/\Cal A_{\Cal P(D,D)}$.}

The endomorphism ring of the object $D$ in the factor category $\Cal C/\Cal A_{\Cal P(D,D)}$ is $\End_R(D)/\Cal A_{\Cal P(D,D)}(D,D)=\End_R(D)/\Cal P(D,D)$, which is a division ring.

\smallskip

{\em Step 2. For every object $A$ of $\Cal C$ with $[A]_{\Cal P}=[D]_{\Cal P}$, the objects $A$ and $D$ are isomorphic objects in the factor category $\Cal C/\Cal A_{\Cal P(D,D)}$.}

Let $A$ be any object of $\Cal C$. Suppose
$[A]_{\Cal P}=[D]_{\Cal P}$. Since $\Cal A_{\Cal P(D,D)}(D,D)=\Cal P(D,D)$, it follows that $\Cal A_{\Cal P(D,D)}\supseteq\Cal P$. In particular, $\Cal A_{\Cal P(D,D)}(A,A)\supseteq\Cal P(A,A)$, which is a maximal ideal. Thus
either $\Cal A_{\Cal P(D,D)}(A,A)=\End_R(A)$ or $\Cal A_{\Cal P(D,D)}(A,A)=\Cal P(A,A).$

In the first case, $1_A\in \Cal A_{\Cal P(D,D)}(A,A)$, so that for every $\alpha\colon D\to A$ and every $\beta\colon A\to~D$ one has that $\beta\alpha\in \Cal P(D,D)$. But $[A]_{\Cal P}=[D]_{\Cal P}$; hence there exist $\alpha\colon D\to A$ and $\beta\colon A\to D$ with $\alpha\notin \Cal P(D,A)$ and $\beta\notin\Cal P(A,D)$. This contradicts the fact that $\Cal P$ is a completely prime ideal.

It shows that
$\Cal A_{\Cal P(D,D)}(A,A)=\Cal P(A,A)$, so that
the endomorphism ring of $A$ in the category
$\Cal C/\Cal A_{\Cal P(D,D)}$ is
$\End_R(A)/\Cal P(A,A)$, which is a division ring. From $[A]_{\Cal P}=[D]_{\Cal P}$, it follows that  $\Cal P(A,D)\ne\Hom(A,D)$ and $\Cal P(D,A)\ne\Hom(D,A)$. Thus there exist $R$-module morphisms $f\colon A\to D$ and $g\colon D\to A$ with $f\notin \Cal P(A,D)$ and $g\notin \Cal P(D,A)$.  But $\Cal P$ is completely prime, so that $fg\notin \Cal P(D,D)$ and $gf\notin\Cal P(A,A)$. Hence $fg+ \Cal P(D,D)$ is an automorphism of $D$ and $gf+\Cal P(A,A)$ is an automorphism of $A$ in the factor category $\Cal \Cal C/\Cal A_{\Cal P(D,D)}$. Thus the morphism $f+ \Cal A_{\Cal P(D,D)}(A,D)$ is both right invertible and left invertible in the category $\Cal C/\Cal A_{\Cal P(D,D)}$. Hence it is an isomorphism in $\Cal C/\Cal A_{\Cal P(D,D)}$. It follows that $A$ and $D$ are isomorphic objects of the category $\Cal C/\Cal A_{\Cal P(D,D)}$.

\smallskip

{\em Step 3. For every object $A$ of $\Cal C$ with $[A]_{\Cal P}\ne[D]_{\Cal P}$, the object $A$ is a zero object in the factor category $\Cal C/\Cal A_{\Cal P(D,D)}$.}

From $[A]_{\Cal P}\ne[D]_{\Cal P}$, it follows that $\Cal P(A,D)=\Hom(A,D)$ or $\Cal P(D,A)={}$ $\Hom(D,A)$. In both cases, for every $\alpha\colon D\to A$ and every $\beta\colon A\to D$, one has that $\beta1_A\alpha\in\Cal P(D,D)$. Thus $1_A\in \Cal A_{\Cal P(D,D)}$.  That is, $A=0$ in $\Cal C/\Cal A_{\Cal P(D,D)}$.

\smallskip

{\em Step 4. Suppose $F:=\{\,i\in I\mid [A]_{\Cal P}=[D]_{\Cal P}\,\}$ is a finite set. Then $P(\prod_{i\in I}A_i)$ is the coproduct in $\Mod R/\Cal A_{\Cal P(D,D)}$ of $|F|$ objects whose endomorphism rings are isomorphic to the division ring $\End_R(D)/\Cal P(D,D)$.}

Since the functor $P$ is additive, one has that $$P(\prod_{i\in I}A_i)\cong P(\oplus_{i\in F}A_i)\coprod P(\prod_{i\in I\setminus F}A_i).$$ By Step 2, $P(\oplus_{i\in F}A_i)$ is the coproduct of $|F|$ objects isomorphic to $D$ in
$$\Mod R/\Cal A_{\Cal P(D,D)}.$$
 Moreover the endomorphism ring of $D$  in $\Mod R/\Cal A_{\Cal P(D,D)}$ is isomorphic to the division ring $\End_R(D)/\Cal P(D,D)$. Let us prove that $$P(\prod_{i\in I\setminus F}A_i)=0$$ in $\Mod R/\Cal A_{\Cal P(D,D)}$; that is, $1_{\prod_{i\in I\setminus F}A_i}\in \Cal A_{\Cal P(D,D)}(\prod_{i\in I\setminus F}A_i,\prod_{i\in I\setminus F}A_i)$. For this, it suffices to show that, for every $\alpha\colon D\to \prod_{i\in I\setminus F}A_i$ and every $\beta\colon \prod_{i\in I\setminus F}A_i\to D$, one has that $\beta\alpha\in {\Cal P(D,D)}$. Now $D$ is slender and $|I|$ is non-measurable, so $\Hom(\prod_{i\in I\setminus F}A_i, D)\cong \bigoplus_{i\in I\setminus F}\Hom(A_i,D)$. Thus if $\pi_j\colon \prod_{i\in I\setminus F}A_i\to A_j$ denotes the canonical projection, then there exists a finite subset $G$ of $I$ disjoint from $F$ and morphisms $\beta_j\colon A_j\to D$ for every $j\in G$ such that $\beta=\sum_{j\in G}\beta_j\pi_j$. As we have seen in Step 3, $A_i$ is a zero object in $\Cal C/\Cal A_{\Cal P(D,D)}$ for every $i\in I\setminus F$. Thus $\beta_j\pi_j\alpha$ is the zero morphism, that is, $\beta\alpha$ is the zero morphism, in $\Mod R/\Cal A_{\Cal P(D,D)}$. We can conclude that $\beta\alpha\in {\Cal P(D,D)}$, as desired.

\smallskip

{\em Step 5. Suppose $\{\,i\in I\mid [A_i]_{\Cal P}=[D]_{\Cal P}\,\}$ is an infinite set. Then, for every integer $n\ge1$, $P(\prod_{i\in I}A_i)$ is the coproduct  of $n$ objects whose endomorphism rings are isomorphic to the division ring $\End_R(D)/\Cal P(D,D)$ and one more object of $\Mod R/\Cal A_{\Cal P(D,D)}$.}

The proof is the same as the first part of the proof of Step 4.

\smallskip

{\em Step 6. Suppose that $\prod_{i\in I}A_i\cong\prod_{j \in J}B_j$. Let $D$ be an object of $\Cal C$. Then the set $\{\,i\in I\mid [A]_{\Cal P}=[D]_{\Cal P}\,\}$ is finite if and only if $\{\,j\in J\mid [B_j]_{\Cal P}=[D]_{\Cal P}\,\}$ is finite. Moreover, in this case, they have the same number of elements.}

Let $\Cal C'$ be any additive category in which idempotents split that contains $\Cal C$, contains $P(\prod_{i\in I\setminus F}A_i)$ for every finite subset $F$ of $I$, and contains $P(\prod_{j \in J\setminus G}B_j)$ for every finite subset $G$ of $J$ \cite[p.~676]{Indiana}. Now apply the Krull-Schmidt-Azumaya Theorem for additive categories \cite[p.~20]{Bass} to the  category $\Cal C'$ and to the object $P(\prod_{i\in I}A_i)\cong P(\prod_{j \in J}B_j)$. Since the endomorphism rings of the non-zero objects $P(A_i)\cong P(B_j)\cong P(D)$ in $\Cal C'$ are division rings, hence local rings, $P(D)^n$ cannot have a direct summand isomorphic to $P(D)^m$ for every $m>n$ by the Krull-Schmidt-Azumaya Theorem. Thus $\{\,i\in I\mid [A]_{\Cal P}=[D]_{\Cal P}\,\}$ is finite if and only if $\{\,j\in J\mid [B_j]_{\Cal P}=[D]_{\Cal P}\,\}$ is finite. Moreover, $P(D)^n\cong P(D)^m$ implies $n=m$. Thus the sets $\{\,i\in I\mid [A]_{\Cal P}=[D]_{\Cal P}\,\}$ and $\{\,j\in J\mid [B_j]_{\Cal P}=[D]_{\Cal P}\,\}$ are equipotent when they are finite.

\smallskip

In order to conclude the proof of the theorem, it is now sufficient to remark that $\{\,i\in I\mid [A]_{\Cal P}=[D]_{\Cal P}\,\}$ and $\{\,j\in J\mid [B_j]_{\Cal P}=[D]_{\Cal P}\,\}$ are always equipotent by hypothesis (b). Now glue the bijections between the sets $\{\,i\in I\mid [A]_{\Cal P}=[D]_{\Cal P}\,\}$ and $\{\,j\in J\mid [B_j]_{\Cal P}=[D]_{\Cal P}\,\}$ to obtain a bijection $\sigma_{\Cal P}$ that preserves the $\Cal P$ classes.
\end{Proof}

As an application of the previous theorem, we get the following Corollary, which is Theorem~2.8 in \cite{Franetic}.

\begin{Cor} {\rm \cite[Theorem~2.8]{Franetic}} Let $R$ be a ring and $\{\,A_i\mid i\in I\,\}$ be a family of slender right $R$-modules with local endomorphism rings.
Let $\{\,B_j\mid j\in J\,\}$ be a family of indecomposable slender right $R$-modules.  Assume that:

{\rm (a)} $|I|$ and $|J|$ are non-measurable cardinals.

{\rm (b)} There are at most countably many mutually isomorphic modules  in each of the two families $\{\,A_i\mid i\in I\,\}$ and $\{\,B_j\mid j\in J\,\}$.

{\rm (c)} The $R$-modules $\prod_{i\in I}A_i$ and $\prod_{j \in J}B_j$ are isomorphic.

Then there exists a bijection $\sigma\colon
I\to J$ such that $A_i\cong B_{\sigma(i)}$ for every $i\in I$. \end{Cor}

\begin{Proof} First of all, let us prove that all the modules $B_j$ also have local endomorphism rings. From (c), each $B_j$ is isomorphic to a direct summand of $\prod_{i\in I}A_i$. As in Step 1 of the proof of Theorem~\ref{nec} \cite[Lemma 1.1]{Franetic}, for each $j\in J$ there exists a finite subset $F_j$ of $I$ such that $B_j$ is a isomorphic to direct summand of $\prod_{i\in F_j}A_i$. Thus $B_j\oplus C_j\cong \prod_{i\in F_j}A_i$, where each $C_j$ is a direct sum of finitely many indecomposable objects, and these indecomposable objects plus $B_j$ are isomorphic, up to a permutation, to the modules $A_i$ with $ i\in F_j$ (Krull-Schmidt-Azumaya Theorem \cite[p.~20]{Bass}). Thus the modules $B_j$ have local endomorphism rings, and now the role of the two families $\{\,A_i\mid i\in I\,\}$ and $\{\,B_j\mid j\in J\,\}$ has become symmetric.

For every pair $A,B$ of objects of $\Cal C$, let $\Cal P(A,B)$ be the set of all the homomorphisms $A\to B$ that are not isomorphisms. It is easily checked that $\Cal P$ is a completely prime ideal of $\Cal C$. Moreover, two objects $A,B$ of $\Cal C$ are isomorphic modules if and only if $[A]_{\Cal P}=[B]_{\Cal P}$. The corollary now follows immediately from Theorem~\ref{buupbp}.\end{Proof}



We conclude the paper with two elementary examples of applications of Theorem~\ref{buupbp}. As a first example, let $R$ be the ring $\Z$ of integers (it could be any other countable principal ideal domain that is not a field). Let\/ $\Cal C$ be the full subcategory of $\Mod \Z$ whose objects are all torsion-free $\Z$-modules $G$ of torsion-free rank $1$ such that $pG\ne G$ for every prime~$p$. There are $2^{\aleph_0}$ pair-wise non-isomorphic such $\Z$-modules (recall that torsion-free abelian groups $G$ of torsion-free rank $1$ are completely determined up to isomorphism by their type $\ty(G)$ \cite[Section~85]{Fuchs2}). All the modules in $\Cal C$ are slender modules \cite[Corollary~III.2.3]{EM}.
For each prime $p$ in $\Z$ and pair $G,H$ of objects of $\Cal C$, set $\Cal P_p(G,H)=p\Hom(G,H)$. The group $\Hom(G,H)$ is torsion-free of rank $1$ and type $\ty(H): \ty(G)$ if $\ty(G)\leq\ty(H)$ \cite[Proposition~85.4]{Fuchs2}, so $\Cal P_p(G,H)<\Hom(G,H)$ in this case. Otherwise, if $\ty(G)\nleq\ty(H)$, then $0=\Cal P_p(G,H)=\Hom(G,H)$. It is then very easy to prove that $\Cal P_p$ is a
 completely prime ideal of $\Cal C$. For every object $A$ of $\Cal C$, $\Cal P_p(A,A)$ is the maximal ideal of $\End_R(A)\cong\Z$ generated by $p$. Thus Theorem~\ref{buupbp} applies to this situation.
Notice that, for every prime $p$ and objects $G,H$ of  $\Cal C$, one has that $[G]_{\Cal P_p} = [H]_{\Cal P_p}$ if and only if $G\cong H$, because if $[G]_{\Cal P_p} = [H]_{\Cal P_p}$, then $\Cal P_p(G,H)<\Hom(G,H)$ and $\Cal P_p(H,G)<\Hom(H,G)$, so that $\ty(G)\leq\ty(H)$ and $\ty(H)\leq\ty(G)$; that is, $\ty(G)=\ty(H)$, and $G\cong H$.

Here is a second example. Recall that a {\em rigid system} of abelian groups is a set $\{\,A_i\mid i\in I\,\}$ of non-zero torsion-free abelian groups for which $\Hom(A_i,A_j)$ is isomorphic to a subgroup of $\Q$ if $i=j$, and is $0$ if $i\ne j$. It is known that there exist rigid systems of abelian groups of finite rank which are homogeneous of type $(0,0,0,\dots)$ \cite[Theorem~88.4]{Fuchs2}. Such groups are torsion-free, reduced and countable, hence slender \cite[Corollary~2.3]{EM}. Let\/ $\Cal C$ be the full subcategory of $\Mod \Z$ with class of objects a rigid system of groups of finite rank homogeneous of type $(0,0,0,\dots)$.
Let $p$ be any prime number and let $\Cal P$ be the completely prime ideal of $\Cal C$ defined, for every $A,B\in\Ob(\Cal C)$, by $\Cal P(A,B)=p\Hom(A,B)$. Then the hypotheses of Theorem~\ref{buupbp} hold. In this case,
$[A]_{\Cal P} = [B]_{\Cal P}$ if and only if $A=B$, for every $A,B\in \Ob(\Cal C)$. Notice that, in these last two examples, the endomorphism rings of the indecomposable direct factors are not local rings.


\end{document}